\documentclass[twoside,11pt]{article}

\usepackage{blindtext}
\usepackage{amsmath}
\usepackage{amssymb}
\usepackage{amsthm}
\usepackage{jmlr2e}

\usepackage{lastpage}
\jmlrheading{23}{2022}{1-\pageref{LastPage}}{1/21; Revised 5/22}{9/22}{21-0000}{Amit Vishwakarma  and  KS Subrahamanian Moosath}

\ShortHeadings{Distance Measure Between GMMs}{}
\firstpageno{1}

\begin{document}

\title{Distance Measure Based on an Embedding of the Manifold of K-Component Gaussian Mixture Models into the Manifold of Symmetric Positive Definite Matrices}

\author{\name Amit Vishwakarma \email amitvishwakarma.22@res.iist.ac.in \\
       \addr Department of Mathematics\\
       Indian Institute of Space Science and Technology\\
       Trivandrum, India
       \AND
       \name KS Subrahamanian Moosath \email smoosath@iist.ac.in \\
       \addr Department of Mathematics\\
       Indian Institute of Space Science and Technology\\
       Trivandrum, India}

\editor{My editor}

\maketitle

\begin{abstract}%
In this paper, a distance between the Gaussian Mixture Models(GMMs) is obtained based on an embedding of the K-component Gaussian Mixture Model into the manifold of the symmetric positive definite matrices. Proof of embedding of K-component GMMs into the manifold of symmetric positive definite matrices is given and shown that it is a submanifold. Then, proved that the manifold of GMMs with the pullback of induced metric is isometric to the submanifold with the induced metric. Through this embedding we obtain a general lower bound for the Fisher-Rao metric. This lower bound is a distance measure on the manifold of GMMs and we employ it for the similarity measure of GMMs.  The effectiveness of this framework is demonstrated through an experiment on standard machine learning benchmarks, achieving accuracy of 98\%, 92\%, and 93.33\% on the UIUC, KTH-TIPS, and UMD texture recognition datasets respectively. 
\end{abstract}


\section{Intordutcion}
 Gaussian Mixture Models (GMMs) are important tools in machine learning, signal processing, and computer vision \citep{rakesh2023moving,8914215} because of their ability to approximate any smooth probability density with a specific non-zero amount of error. Although GMMs are widely used, a key challenge lies in the lack of closed-form expressions for similarity measures and geodesic between them, which impacts the effectiveness of GMM-based applications, particularly when comparing and analyzing them in higher dimensional spaces.

For univariate and multivariate normal distributions, several well-established methods to measure similarities already exists. Particularly in the multivariate case, there exists a closed-form expression for the Kullback-Leibler divergence between two normal distributions \citep{pardo2018statistical}, providing a natural way to compare them. Calvo and Oller obtained a robust geometric framework by embedding multivariate normal distributions $N(x;\mu,\Sigma)$ into the manifold of symmetric positive definite matrices $SPD_{(n+1)}(\mathbb{R})$ and  achieved a lower bound for the Rao distance which itself is a metric and also derived general statistical test of hypothesis \citep{CALVO1990223}. 

The Kullback-Leibler (KL) divergence \citep{kullback1951information}, a fundamental information-theoretic measure, is often the preferred choice to analyze the distributional differences. However, lack of a closed-form expression for the KL divergence between GMMs has led researchers to explore various alternatives. These include developing KL lower and upper bounds \citep{4218101}, approximation techniques (\cite{peter2006shape}), Monte Carlo sampling methods and closed-form Cauchy-Schwarz divergence for a mixture of Gaussians \citep{Kampa2011ClosedformCP}. However, these methods often face challenges with computational efficiency and accuracy, particularly when dealing with GMMs that have multiple low-probability components or significant overlapping regions.

Recent works have attempted to address the computational challenges in
GMM similarity measures through various approaches. In \cite{li2013novel}, a novel Earth Mover's Distance (EMD) methodology for GMM matching is introduced, proposing a sparse representation-based EMD (SR-EMD) that exploits the inherent sparsity of the transportation problem. They also developed ground distances between components of GMMs based on information geometry, viewing the space of Gaussians as either embedded in a Lie group or as a product of Lie groups. While this approach improved efficiency and robustness compared to conventional EMD, it still depends on component-wise analysis and did not fully capture the geometric structure of the entire GMM.

\cite{popovic2021measure} explored the embedding of components of GMMs into symmetric positive definite matrices while preserving KL divergence properties. Their work focused primarily on computational efficiency. This method also not utilizing the full geometric structure of the GMM manifold. More recently, \cite{popovic2023measure} proposed a geometry-aware dimensionality reduction technique for GMMs that preserves local neighborhood information through Distance Preservation to Local Mean (DPLM). This method achieved significant computational efficiency by operating in a reduced dimensional space. However, a fundamental limitation in \cite{popovic2023measure}, \cite{popovic2021measure} is handling GMMs through component-wise analysis and embedding the multivariate Gaussian components $N(x;\mu_k,\Sigma_k)$ of GMMs into the manifold of symmetric positive definite matrices $SPD_{(n+1)}(\mathbb{R})$, potentially losing important statistical relationships between components and fails to use the geometric structure of the manifold of GMMs.

This paper addresses these limitations by embedding the manifold of K-component GMMs into the manifold of symmetric positive definite matrices, instead of analyzing individual components. Our approach maps K-component GMMs into a higher-dimensional manifold of symmetric positive definite matrices while preserving their intrinsic probabilistic structure. It is proved that the manifold of K-component GMMs forms a submanifold of the manifold of symmetric positive definite matrices endowed with a natural Riemannian metric induced from the ambient manifold. We also establish several key theoretical results: (1) the embedded manifold is a submanifold of dimension $ \frac{K}{2}(n+1)(n+2)-1$, (2) the manifold of GMM is isometric to its image under this embedding, where the isometry is defined with respect to the affine-invariant Riemannian metric, (3) while generally non-geodesic, the embedded submanifold becomes geodesic under conditions of fixed means and uniform mixing coefficients, and (4) the Fisher-Rao distance in the manifold of GMMs is always greater than or equal to the affine metric in the embedded space, with equality holding under the aforementioned conditions.

This approach provides several key advantages over existing approaches. Unike SR-EMD \citep{li2013novel} and DPLM \citep{popovic2023measure} this method takes care of the geometry of the K-component GMMs rather than the component-wise approach, which helps to a more comprehensive understanding of the manifold of GMMs. Furthermore, while the DPLM approach of \cite{popovic2023measure} preserves the local neighborhood information, this framework maintains both local and global geometric properties of the GMM manifold, providing a complete understanding of the geometry of the manifold of GMMs and enabling more robust applications in pattern recognition, texture recognition and statistical inference.

The effectiveness of our geometric framework is demonstrated through experiments on standard texture recognition benchmarks. Our method achieves significant improvements in classification accuracy 98.82\%, 92\%, and 93.33\% on UIUC (the University of Illinois Urbana-Champaign), KTH-TIPS ( the Royal Institute of Technology - Textures under varying Illumination, Pose and Scale), and UMD ( University of Maryland texture datasets) respectively. This validate both the theoretical and practical use of this framework for real-world machine-learning tasks.

\section{Related Works}

In this section, we discussed various approaches in the similairty measure of GMMs. We first discuss classical information-theoretic approaches based on KL divergence and their limitations. We then examine various approximation techniques that have been proposed to address these limitations, followed by recent geometric approaches that form the basis for our work.\\

The Kullback-Leibler (KL) divergence serves as a fundamental measure of dissimilarity between probability distributions in information theory and statistical learning. For two probability distributions p and q, the KL divergence is defined as,
\begin{equation}
    KL(p||q) = \int_{\mathbb{R}^d} p(x) \log \frac{p(x)}{q(x)} dx.
\end{equation}
While this measure has a closed-form expression for individual Gaussian distributions, given by,

\begin{equation}
    KL(\hat{p}\|\hat{q}) = \frac{1}{2}\left[\log \frac{|\Sigma_{\hat{q}}|}{|\Sigma_{\hat{p}}|} + \text{Tr}[\Sigma_{\hat{q}}^{-1}\Sigma_{\hat{p}}] + (\mu_{\hat{p}} - \mu_{\hat{q}})^T\Sigma_{\hat{q}}^{-1}(\mu_{\hat{p}} - \mu_{\hat{q}}) - d\right],
\end{equation}
no such analytical solution exists for GMMs.
The most direct approach to estimating KL divergence between GMMs is through Monte Carlo sampling \citep{4218101}. The idea is for GMMs $p = \sum_{i=1}^n \alpha_i p_i$ and $q = \sum_{j=1}^m \beta_j q_j$, using N i.i.d. samples drawn from p(x) the divergence is,

\begin{equation}
    KL_{MC}(p\|q) = \frac{1}{N}\sum_{i=1}^N \log \frac{p(x_i)}{q(x_i)}.
\end{equation}

While theoretically accurate as N approaches infinity, this method is computationally expensive for real-world applications, particularly when high accuracy is required.\\

Several component-wise approximations have been proposed to address the computational challenges. The Weighted Average approximation \citep{1238387} leverages the convexity of KL divergence to obtain
\begin{equation}
    KL_{WA}(p\|q) \approx \sum_{i,j} \alpha_i \beta_j KL(p_i\|\|q_j).
\end{equation}
While computationally efficient, this approximation becomes inaccurate when mixture components are well-separated. \\
The Matching-Based \citep{1238387} approximation improves upon this by assuming dominant contributions from the closest components
\begin{equation}
    KL_{MB}(p\|q) \approx \sum_i \alpha_i \min_j\left[KL(p_i\|\|q_j) + \log\left(\frac{\alpha_i}{\beta_j}\right)\right].
\end{equation}
This method performs well when components are distant but struggles with significant component overlap.\\
The Variational Approximation \citep{Durrieu2012LowerAU} offers a more sophisticated approach through
\begin{equation}
KL_{VA}(p||q) = \sum_i \alpha_i \frac{\sum_{i'} \alpha_{i'} e^{-KL(p_i||p_{i'})}}{\sum_j \beta_j e^{-KL(p_i||q_j)}}.
\end{equation}
This formulation better handles component overlap by considering relationships between all Gaussian components in both mixtures. However, it still relies on component-wise analysis which may not fully capture the geometric structure of the GMM manifold.\\

More recent approaches have investigated nonlinear dimensionality reduction techniques to improve computational efficiency while preserving geometric structure. Notably,  \cite{popovic2023measure} proposed a geometry-aware dimensionality reduction framework based on Distance Preservation to the Local Mean (DPLM) for the manifold of symmetric positive definite matrices. This approach embeds individual Gaussian components $N(x;\mu_k,\Sigma_k)$ into $SPD_{(n+1)}(\mathbb{R})$ through the mapping,

\begin{equation}
g \hookrightarrow P = |\Sigma|^{-\frac{1}{d+1}} 
\begin{bmatrix} 
\Sigma + \mu\mu^T & \mu \\ 
\mu^T & 1 
\end{bmatrix}
\end{equation}

The DPLM algorithm then finds a projection matrix W that maps these embeddings into a lower-dimensional space $SPD_{(l)}$ where $l < d$ while preserving local geometric structure. This is achieved by minimizing
$H(W) = \sum_{i=1}^N \sum_{k=1}^K |\delta_{ld}^2(P_{i,k}, \hat{N_i}) - \delta_{ld}^2(W^TP_{i,k}W, W^T\hat{N_i}W)|$
subject to $W^T \times W = I_l$, where $\delta_{ld}^2$ is a distance measure derived from the Jensen-Bregman LogDet divergence:
$\delta_{ld}^2(P,Q) = \sqrt{J(P,Q)}$ with $J(P,Q) = \text{logdet}(\frac{P + Q}{2}) - \frac{1}{2}\text{logdet}(PQ)$
and $\hat{N_i}$ represents the Riemannian mean of the K-nearest neighbors of each embedding $P_i$. While this approach achieves significant computational gains through dimensionality reduction and preserves local geometric relationships, it shares a fundamental limitation with other existing methods - it analyzes GMMs through component-wise embeddings rather than treating them as unified probabilistic objects. This component-wise treatment potentially loses important statistical relationships between mixture components and fails to leverage the rich geometric structure of the GMM manifold. Also, the choice of K nearest neighbors and projection dimension $l$ becomes critical for performance, requiring careful tuning for each application.
The limitations of these approaches motivate our development of a unified geometric framework that embeds the K-component GMMs into the manifold of symmetric positive definite matrices while preserving their complete probabilistic structure. This allows us to explore the geometry of the embedded space for efficient computation of similarity measures while maintaining the statistical relationships between mixture components.

\section{Preliminaries}
In this section, the fundamental concepts of statistical manifold and the geometric structure are introduced, including the Fisher information metric and divergence measure. The stastical manifold structure of GMMs and the manifold of symmetric positive definite matrices are also given.

\subsection{Statistical manifold}
Now, we describe the manifold structure and its geometry for a statistical model. 

Let $(\mathcal{X}, \Sigma, p)$ be a probability space, where $\mathcal{X} \subseteq \mathbb{R}^{n}$. Consider a family $\mathcal{S}$ of probability distributions on $\mathcal{X}$. Suppose each element of $\mathcal{S}$ can be parametrized using $n$ real-valued variables $(\theta^{1},...,\theta^{n})$ so that 
\begin{equation}
\mathcal{S}=\lbrace p_{\theta}=p(x;\theta)\;/ \;\theta=(\theta^{1},...,\theta^{n}) \in \Theta \rbrace
\end{equation}
where $\Theta$ is an open subset of $\mathbb{R}^{n}$ and the mapping $\theta \mapsto p_{\theta}$ is injective. The family $\mathcal{S}$ is called an $n$-dimensional \textbf{statistical model} or a\textbf{ parametric model}. We often write as $\mathcal{S}=\lbrace p_{\theta}\rbrace.$

For a model $\mathcal{S}=\lbrace p_{\theta} \;/ \; \theta \in \Theta \rbrace$, the mapping $\varphi:\mathcal{S}\longrightarrow \mathbb{R}^{n}$ defined by $\varphi(p_{\theta})=\theta$ allows us to consider $\varphi=(\theta^{i})$ as a coordinate system for $\mathcal{S}$. Suppose there is a $\mathit{c^{\infty}}$ diffeomorphism $\psi:\Theta\longrightarrow \psi(\Theta)$, where $\psi(\Theta)$ is an open subset of $\mathbb{R}^{n}$. Then, if we use $\rho=\psi(\theta)$ instead of $\theta$ as our parameter, we obtain $\mathcal{S}=\lbrace p_{\psi^{-1}(\rho)}\;\vert \; \rho \in \psi(\Theta) \rbrace$. This expresses the same family of probability distributions $\mathcal{S}=\lbrace p_{\theta} \rbrace$. If parametrizations which are $\mathit{c^{\infty}}$ diffeomorphic to each other is considered to be equivalent then $\mathcal{S}$ is a $\mathit{c^{\infty}}$ differentiable manifold, called the \textbf{statistical manifold}.
The tangent space $T_{\theta}(\mathcal{S})$  to the statistical manifold  $ \mathcal{S} $ at a point $p_{\theta}$ is spanned by 
\begin{equation}
\lbrace \frac{\partial}{\partial \theta^{i}} \vert_{p_{\theta}} ; i=1,...,n\rbrace.
\end{equation}
The score functions $ \lbrace \frac{\partial}{\partial \theta_{i}} \log p(x;\theta) ;i=1,...,n\rbrace$ are assumed to be linearly independent functions in $x$. Denote $\frac{\partial}{\partial \theta_{i}} \log p(x; \theta)$
 by $ \partial_{i} \ell(x, \theta) $.
The vector space $T_{\theta}^{1}(S)$ spanned by $\{\partial_{i}l(x,\theta): i=1,..n\}$ is isomorphic with $T_{\theta}(S).$ Note that $T_{\theta}(S)$ is the differential operator representation of the tangent space  and $T_{\theta}^{1}(S)$ the random variable representation of the tangent space.
\\
\noindent The expectation $E_{\theta}[\partial_{i} \ell ({x;\theta})] = \int\partial_{i} \ell({x;\theta}) p(x;\theta)dx=0,$ since $p(x;\theta)$ satisfies $\int p(x;\theta)dx=1$ and so $E_{\theta}[A(x)]=0$ for  $A(x) \in T_{\theta}^{1}(\mathcal{S})$.
This expectation induces an inner product on $\mathcal{S}$ in a  natural way as $<A,B>_{\theta}\;=\;E_{\theta}[A(x)B(x)].$
Denote $\frac{\partial}{\partial \theta_{i}}$ by $\partial_{i}$
the inner product of the basis vectors $\partial_{i}$ and $\partial_{j}$ is 
\begin{eqnarray}
g_{ij}(\theta) &=& \langle \partial_{i},\partial_{j} \rangle_{\theta} \nonumber \\
&=& E_{\theta}[\partial_{i} \ell (x;\theta) \partial_{j} \ell (x;\theta)] \nonumber \\
&=& \int \partial_{i} \ell(x;\theta) \partial_{j} \ell(x;\theta) p(x;\theta) dx \label{eq:k1}.
\end{eqnarray}
 Note that the matrix $G(\theta)=(g_{ij}(\theta))$ is symmetric. For the vector $c=[c^{1},...,c^{n}]^{t}$\\
$c^{t}G(\theta)c=\int \lbrace \sum_{i=1}^{n} c^{i} \partial_{i}\ell(x;\theta)\rbrace ^2 p(x;\theta) dx \geq 0.$ and since $\lbrace \partial_{1}\ell(x;\theta),...,\partial_{n} \ell(x;\theta) \rbrace $  are linearly independent, $G$ is positive definite. 
Hence $g=<,>$ defined above is a Riemannian metric on the statistical manifold $\mathcal{S}$,  called the \textbf{Fisher information metric or Fisher-Rao metric} \citep{amari2000methods}.

\subsection{Divergence Measures}
On a statistical manifold, divergences are use to measure dissimilarity between probability distributions. Divergence is a distance-like measure between two points (probability density functions) on a statistical manifold. The divergence $D$ on $S$ is defined as $D(.||.):S \times S \to \mathbb{R}$, a smooth function satisfying, for any $p,q\in S$
$$D(p||q) \geq 0\text{ and }D(p||q)=0\text{ iff }p=q.$$

Among various divergence measures, the Kullback-Leibler (KL) divergence plays a central role in information geometry. The KL divergence is defined as \citep{kullback1951information},
\begin{equation}
    D_{KL}(p || q) = \int p(x; \theta_1) \log\frac{p(x; \theta_1)}{q(x; \theta_2)} dx.
\end{equation}
This measure, while not symmetric and not satisfying the triangle inequality, has deep connections to the Fisher information metric and provides a fundamental tool for comparing probability distributions.

\subsection{Gaussian Mixture Models}
Gaussian Mixture Models is a powerful tool in statistical modeling, offering a flexible approach to representing complex, multidimensional probability distributions. These models are widely applied across different fields, such as machine learning, pattern recognition, and signal processing etc \citep{4218101}.

\begin{definition}
 A Gaussian Mixture Model is a probabilistic model that represents a distribution as a weighted sum of Gaussian component densities. For a $n$-dimensional random vector $\mathbf{x}$, a GMM with $K$ components is defined by the probability density function
\begin{equation}
    p(\mathbf{x}) = \sum_{k=1}^{K} \pi_k \mathcal{N}(\mathbf{x}; \mu_k, \Sigma_k)
\end{equation}
where $\pi_k$ are the mixture weights satisfying $\sum_{k=1}^{K} \pi_k = 1$ and $\pi_k \geq 0$ for all $k$, and $\mathcal{N}(\mathbf{x}; \mu_k, \Sigma_k)$ is the probability density function of a multivariate Gaussian distribution with mean vector $\mu_k$ and covariance matrix $\Sigma_k$.
\end{definition}

The flexibility of GMMs arises from their ability to approximate any continuous density function with a high accuracy.

\begin{theorem}[Universal Approximation]
Let $f$ be a continuous probability density function on $\mathbb{R}^d$. For any $\epsilon > 0$, there exists a Gaussian Mixture Model $p(\mathbf{x})$ with $K$ components such that
.\begin{equation}
    \sup_{\mathbf{x} \in \mathbb{R}^d} \big|f(\mathbf{x}) - p(\mathbf{x})\big| < \epsilon
    \tag*{\citep{bengio2017deep}}
\end{equation}

\end{theorem}

This approximation capability makes GMMs a powerful tool in machine learning and statistical modeling. To fully characterize a GMM, we need to specify its parameters. For an $n$-dimensional GMM with $K$ components, the parameter space $\Theta$ can be represented as
\begin{equation}
    {\Theta} = \{\theta\in \mathbb{R}^{ \frac{K}{2}(n+1)(n+2)-1}:\theta=(\pi_k,\mu_k,\Sigma_k)_{k=1}^{K}\}
\end{equation}
where,
 $\pi_k$ is the first $K-1$ mixture weights (the $K^{th}$ weight is determined by the constraint $\sum_{k=1}^K \pi_k = 1$), $\mu_{k}$ is the mean vector, $\Sigma=[\sigma_{ij}]$ is the $n\times n$ covariance matrix.

Let $\mathcal{M}=\{p(x;\theta):\theta\in \Theta\},$ where $p(\mathbf{x}; \theta)$ is the probability density function of the GMMs with parameters $\theta\in \Theta$ be the set of all GMMs with K components in $n$-dimensional space. This statistical model $\mathcal{M}$ forms a statistical manifold of dimension
\begin{equation}
    \frac{K}{2}(n+1)(n+2)-1.
\end{equation}

 Each point on this manifold corresponds to a specific GMM with a unique set of parameters $\theta.$ So, in the statistical manifold $\mathcal{M}$ each point $p(x;\theta)$ is represented by the parameter set $\theta\in \Theta$.

Understanding the geometric properties of the GMM manifold is crucial for developing meaningful distance measures and statistical inference procedures. The main challenge in embedding the statistical manifold of GMMs into the manifold of symmetric positive definite matrices lies in its parameterization. For a single Gaussian distribution, the embedding is straightforward, defined by its mean and covariance matrix \citep{CALVO1990223}. We extend the Calvo-Oller embedding to the more complex case of Gaussian Mixture Models and try to understand the geometry of GMM manifold.

\subsection{ Manifold of Symmetric Positive Definite Matrices}

A symmetric real $n \times n$ matrix $A$ is said to be positive definite if $x^T A x > 0$ for all non-zero $x \in \mathbb{R}^n$. The set of all $n \times n$ symmetric positive definite matrices is denoted by $SPD_n(\mathbb{R})$. The space $SPD_n(\mathbb{R})$ is an open subset of the space of $n \times n$ symmetric matrices $S_n(\mathbb{R})$, and thus inherits a smooth manifold structure of dimension $n(n+1)/2$. It is also a subset of the space of positive definite matrices, which is denoted by $P_n(\mathbb{R})$.
The space of symmetric matrices can be considered as the minimum linear extension of the symmetric positive definite manifold. So the tangent space at a point $A \in SPD_n(\mathbb{R})$ is identified with the space of symmetric matrices
 \[
T_A SPD_n(\mathbb{R}) = \{ X \in \mathbb{R}^{n \times n} : X = X^T \}.
\]

The geometry of $SPD_n(\mathbb{R})$ can be described using the differential metric,
\[
ds^2 = \frac{1}{2}\text{Tr}\{(S^{-1} dS)^2\}.
\]
This metric provides the infinitesimal squared distance on the manifold, where $S \in SPD_n(\mathbb{R})$ and $ds$ is an infinitesimal perturbation in the tangent space.

The manifold $SPD_n(\mathbb{R})$ can be equipped with several natural Riemannian metrics. One important metric is the affine-invariant metric.

\begin{definition}[Affine-Invariant Metric]
For $A \in SPD_n(\mathbb{R})$ and $X, Y \in T_A SPD_n(\mathbb{R})$, the affine-invariant metric is defined as
\[
\rho_A(X, Y) = \text{tr}(A^{-1} X A^{-1} Y).
\]
This metric is invariant under the action of $GL_n(\mathbb{R})$ on $SPD_n(\mathbb{R})$ given by $A \mapsto G^T A G$ for $G \in GL_n(\mathbb{R})$ \citep{Burbea1984InformativeGO}.
\end{definition}

The geodesics and geodesic distance under this metric have closed-form expressions. The geodesic between $A, B \in SPD_n(\mathbb{R})$ is given by,
    \[
    \gamma(t) = A^{1/2} (A^{-1/2} B A^{-1/2})^t A^{1/2}, \quad t \in [0, 1]
    \] The geodesic distance between $A$ and $B$ is given by,
    \[
    d(A, B) = \| \log(A^{-1/2} B A^{-1/2}) \|_F,
    \]
    where $\log$ is the matrix logarithm and $\| \cdot \|_F$ is the Frobenius norm.
\
The exponential and logarithmic maps play crucial role in the geometry of $SPD_n(\mathbb{R})$. For $A \in SPD_n(\mathbb{R})$ and $X \in T_A SPD_n(\mathbb{R}),$ the exponential map is given by $
    \text{Exp}_A(X) = A^{1/2}\exp(A^{-1/2} X A^{-1/2})A^{1/2}
    .$
     The logarithmic map from $A$ to $B \in SPD_n(\mathbb{R})$ is $
    \text{Log}_A(B) = A^{1/2}\log(A^{-1/2} B P^{-1/2})A^{1/2}$ where, $\exp$ and $\log$ denote the matrix exponential and logarithm, respectively.

\section{ Embedding of the statistical manifold $\mathcal{M}$ of GMMs into the manifold $SPD_{K(n+1)}(\mathbb{R})$}
In this section, we prove the embedding of the statistical manifold $\mathcal{M}$ of K-component GMMs into the manifold  $SPD_{K(n+1)}(\mathbb{R})$ of symmetric positive definite matrices. First, the map form $\mathcal{M}$ into $SPD_{K(n+1)}(\mathbb{R})$ is defined in theorem 4.

\begin{theorem}\label{app:S}
The map $f: \mathcal{M}\to SPD_{K(n+1)}(R)$ by $f(\theta)=S$ where $\theta=(\pi_k,\mu_k,\Sigma_k)$ and the matrix $S$ is given by

\begin{equation}
S = \begin{bmatrix}
A & X \\
X^T & B
\end{bmatrix},
\label{eq:distance}
\end{equation}
where
\[
A = \begin{bmatrix}
\Sigma_1 + \pi_1 \mu_1 \mu_1^T & 0 & \cdots & 0 \\
0 & \Sigma_2 + \pi_2 \mu_2 \mu_2^T & \cdots & 0 \\
\vdots & \vdots & \ddots & \vdots \\
0 & 0 & \cdots & \Sigma_k + \pi_k \mu_k \mu_k^T
\end{bmatrix}
\]
\[
X = \begin{bmatrix}
\pi_1 \mu_1 & 0 & \cdots & 0 \\
0 & \pi_2 \mu_2 & \cdots & 0 \\
\vdots & \vdots & \ddots & \vdots \\
0 & 0 & \cdots & \pi_k \mu_k
\end{bmatrix}
\]
\[
B = \begin{bmatrix}
\pi_1  & 0 & \cdots & 0 \\
0 & \pi_2  & \cdots & 0 \\
\vdots & \vdots & \ddots & \vdots \\
0 & 0 & \cdots & \pi_k 
\end{bmatrix}
\]
is well defined, where   $\Sigma_k$ are $n \times n$ positive definite covariance matrices, $\mu_k$ are $n \times 1$ mean vectors, and $\pi_k$ are mixing coefficients, for $k = 1, \ldots, K$.
\end{theorem}

    Proof is given in Appendix \ref{app:SPDT}.\\
Next we give the differential metric expression for such matrices.
\begin{theorem}
For any $S \in SPD_{K(n+1)}(\mathbb{R})$ of the form
\[
S = \begin{pmatrix}
A & X \\
X^T & B
\end{pmatrix},
\]
where $A,X,X^{T},B$ are as given in theorem 4.
 the differential metric is given by,
\[
ds^2 = \sum_{k=1}^K \left[ \frac{1}{2} \left(\frac{d\pi_K}{\pi_K}\right)^2 + \pi_K d\mu_k^T \Sigma_k^{-1} d\mu_k + \frac{1}{2} \text{tr} \{ (\Sigma_k^{-1} d\Sigma_k)^2 \} \right].
\]
\end{theorem}
    Proof is given in Appendix \ref{app:differential}.\\

Now, we discuss the embedding of the statistical manifold $\mathcal{M}$ of K-component GMMs into the manifold $SPD_{K(n+1)}(\mathbb{R})$.

\begin{theorem}
The map $f:\mathcal{M}\to SPD_{K(n+1)}(\mathbb{R})$ given by $f(\theta)=\begin{pmatrix}
A & X \\
X^T & B
\end{pmatrix} $ satistfies
\begin{enumerate}
\item $f$ is an embedding and $f(\mathcal{M})$ is a submanifold of $SPD_{K(n+1)}(\mathbb{R})$ of dimension $ \frac{K}{2}(n+1)(n+2)-1$.
\item The manifold $\mathcal{M}$ with the pullback of the induced metric is isometric to its image $f(\mathcal{M})$ with the induced metric.
\item In general, \( f(\mathcal{M}) \) is a non-geodesic submanifold of \( SPD_{K(n+1)}(\mathbb{R}) \).
\item $f(\mathcal{M})$ becomes the geodesic submanifold of \( SPD_{K(n+1)}(\mathbb{R}) \) when the mean vectors \( \mu_k \) are fixed and the mixing coefficients \( \pi_k \) are uniform and constant.

\end{enumerate}
\end{theorem}

    Proof is given in Appendix \ref{app: Geometric}.\\

\begin{theorem}
Let $(\mathcal{M},g_{\mathcal{M}})$ be the statistical manifold of GMMs with Fisher-Rao metric and $(SPD_{K(n+1)}(\mathbb{R}),\rho)$ be the manifold of symmetric positive definite matrices with affine-invariant metric $\rho$. Let $\rho_{f(\mathcal{M})}$ be the induced affine-invariant metric on the submanifold $f(\mathcal{M})$. Then $g_{\mathcal{M}} \geq \rho_{f(\mathcal{M})}$. The equality holds when mean vectors are fixed and mixing coefficients are uniform.
\end{theorem}
Proof is given in Appendix \ref{app:inequality}.\\

\textbf{Note:} The theoretical framework presented above allows use of the affine-invariant metric for measuring distances between Gaussian Mixture Models. Manifold of GMMs is embedded in $SPD_{K(n+1)}(\mathbb{R})$ and is a submanifold of dimension $\frac{K}{2}(n+1)(n+2)-1$. $(\mathcal{M},f^{*}\rho_{f(\mathcal{M})})$ is isometric to $(f(\mathcal{M}),\rho_{f(\mathcal{M})})$ and also obtained a general lower bound for the Fisher-Rao metric. This lower bound is a distance measure on the manifold of GMMs and we employ it for the simialrity measure of GMMs.

\section{Result and Discussion}
In this section, we do the texture classification using the metric obtained on the GMM manifold as discussed in the above framework. The experiment is done using three widely-used texture recognition benchmarks: the UIUC (University of Illinois Urbana-Champaign) dataset \citep{1453514}, KTH-TIPS (Textures under varying Illumination, Pose and Scale) dataset \citep{KTHTipsDatabase}, and UMD (University of Maryland) dataset \citep{Xu2009ViewpointInvariant}. Our evaluation compares the proposed method against established KL-based GMM similarity measures, specifically the Weighted Average (\text{KL\_WA}) and Matching-Based (\text{KL\_MB}) approximations, which are proven for reliability across various pattern recognition tasks.

\begin{figure}[h!]
    \centering
    \includegraphics[width=\linewidth]{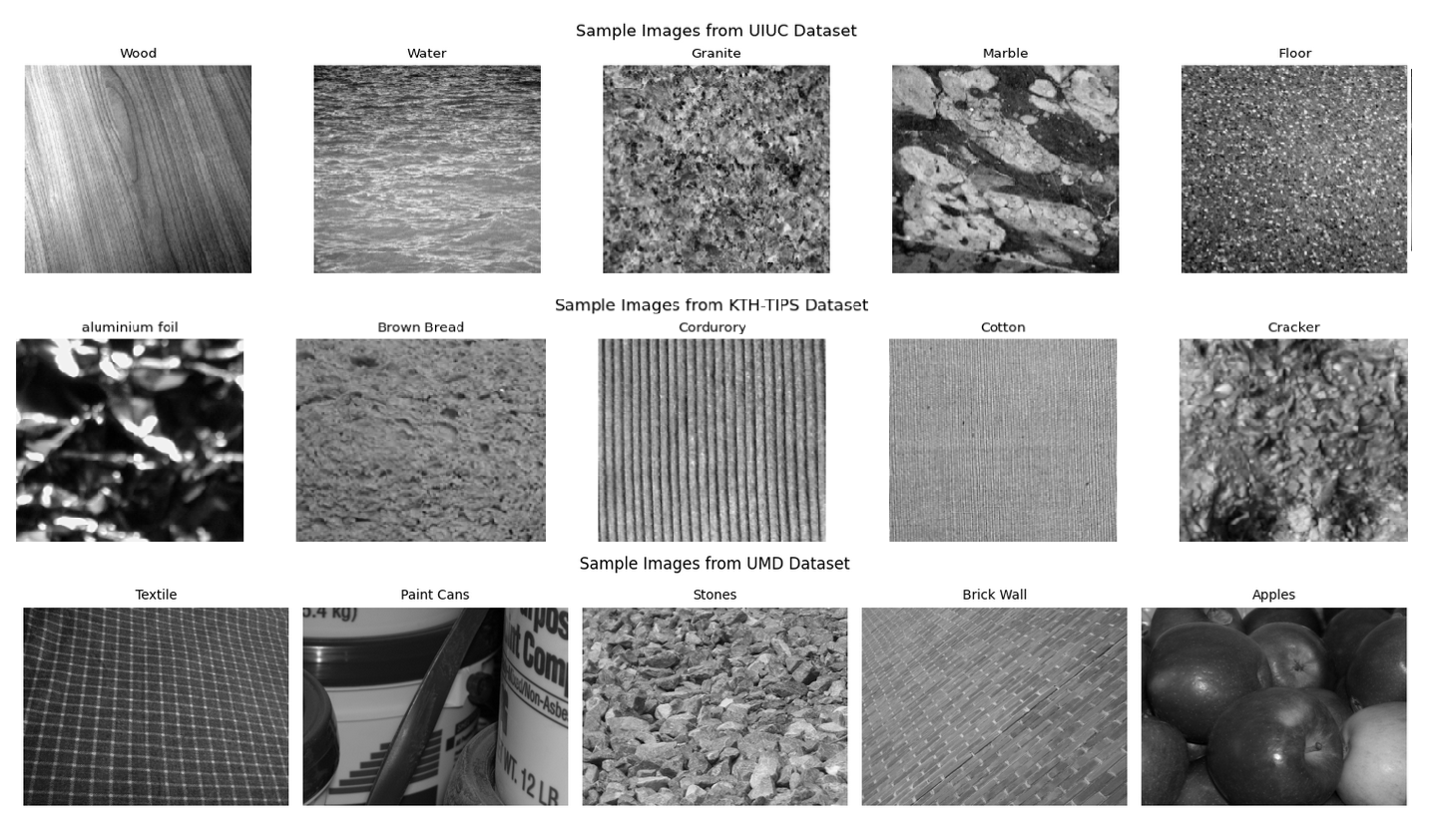} 
    \caption{Sample data of the UIUC, KTH-TIPS, and UMD datasets.}
    \label{fig:sample_data}
\end{figure}

\subsection{Dataset Description and Preprocessing}
Our experiments were conducted using three standard texture databases that are widely used in computer vision and pattern recognition research. The UIUC texture database provided five distinctive texture classes (wood, water, granite, marble, and floor), with images captured at high resolution ($640\times 480$ pixels) to ensure detailed texture information. Each class contains multiple discription of the samples, enabling robust evaluation of our method's invariance properties.
From the KTH-TIPS database, which is specifically designed to evaluate texture recognition under varying illumination, pose, and scale conditions, we selected five representative texture classes (aluminum foil, brown bread, corduroy, cotton, and cracker). The images were standardized to $200\times 200$ pixels to ensure uniform processing while preserving essential texture characteristics. This dataset is particularly challenging due to its systematic variations in imaging conditions, providing a rigorous test of our method's robustness.
The UMD texture dataset contributed five diverse categories: paint cans (class 2), stones (class 3), brick walls (class 8), apples (class 9), and textile patterns (class 12). These images, sampled at $1280\times 960$ pixels, offer a wide range of natural and man-made textures with significant intra-class variations. The high resolution and diverse texture patterns in this dataset make it particularly suitable for evaluating our unified GMM embedding framework's ability to capture complex texture characteristics.

\subsection{Feature Extraction Process}

Region covariance descriptors \citep{Sivalingam2010TensorSparseCoding} were employed as texture features due to their demonstrated effectiveness in texture recognition tasks. For each texture image, the feature extraction process consisted of multiple stages. First, images were split into four equal-sized quadrants to capture local texture variations. For each quadrant, patches were extracted with dimensions specific to each dataset: $128\times 128$ pixels (step 16) for UIUC, $40\times 40$ pixels (step 5) for KTH-TIPS, and $256\times 256$ pixels (step 32) for UMD.

For every pixel position $(x,y)$ within each patch, a five-dimensional feature vector $f(x,y)$ was computed:

\[
f(x,y) = [I(x,y), |I_x(x,y)|, |I_y(x,y)|, |I_{xx}(x,y)|, |I_{yy}(x,y)|]
\]

where $I(x,y)$ represents the pixel intensity, $I_x$ and $I_y$ denote the first-order derivatives in $x$ and $y$ directions computed using Sobel operators, and $I_{xx}$ and $I_{yy}$ represent the second-order derivatives. The absolute values of derivatives were used to ensure robustness to illumination changes.

For each patch, given a set of $m$ feature vectors $\{f_i\}_{i=1}^m$, a $5\times 5$ covariance matrix $C$ was computed:

\[
C = \frac{1}{m-1}\sum_{i=1}^m (f_i - \mu)(f_i - \mu)^T
\]

where $\mu$ is the mean vector of the features. Due to the symmetry of covariance matrices, only the upper triangular part was retained, resulting in a $d$-dimensional feature vector where $d = \frac{n(n+1)}{2} = 15$ for $n=5$ features. This vectorization process transforms each patch into a compact representation capturing the correlations between different texture features.

\subsection{GMM Parameter Estimation}
The parameters of Gaussian Mixture Models were estimated using the Expectation Maximization (EM) algorithm \citep{Webb2003StatisticalPatternRecognition} applied over the set of extracted feature vectors. Each quadrant was represented by a GMM with five components using full covariance matrices. Given a set of feature vectors extracted from patches within a quadrant, each GMM was parameterized by mixture weights $\pi_k$, mean vectors $\mu_k$, and covariance matrices $\Sigma_k$, where $k = {1,\dots,5}$ is the number of components. The mixture weights satisfy the constraint $\sum_{k=1}^5 \pi_k = 1$ with $\pi_k \geq 0$, and each component follows a normal distribution $\mathcal{N}(\mu_k, \Sigma_k)$ in the 15-dimensional feature space.

\subsection{Embedding of GMMs}
Each GMM is embedded into the space of symmetric positive definite matrices by constructing a block matrix $S$ as defined in equation \ref{eq:distance} that preserves both individual component characteristics and their inter-relationships. 

This embedding ensures that the resulting matrix $S$ is symmetric positive definite and preserves the complete probabilistic and geometric structure of the original GMM. Each quadrant's GMM is independently embedded into this symmetric positive definite matrix representation, maintaining the hierarchical structure of the texture analysis from quadrant-level to image-level comparisons.

\subsection{Classification Framework}
The classification of texture images is based on the affine-invariant metric on the manifold of symmetric positive definite matrices. For two SPD matrices P and Q, this metric is defined as $d(P,Q) = \|log(P^{-1/2} Q P^{-1/2})\|_F$, where $\log$ denotes the matrix logarithm and $\|\cdot\|_F$ is the Frobenius norm. To optimize computation, we employ two key strategies: (1) a precomputation phase where $P^{-1/2}$ is computed via eigendecomposition for each training SPD matrix, with a small regularization term $\epsilon$ added for numerical stability, and (2) batch processing of test data to enable efficient parallel computation. During distance computation, each matrix pair undergoes regularization followed by symmetrization of the intermediate matrix $M = P^{-1/2} Q P^{-1/2}$ to ensure robust eigenvalue computation.

The classification follows a two step voting scheme where each quadrant's SPD matrix of the test data is first classified using K-nearest neighbors (K=5). The distances to all training SPD matrices are computed using the affine-invariant metric, and the label is assigned through majority voting among the K nearest neighbors. The final image classification is then determined by aggregating these quadrant-level predictions through a second majority voting step across the four quadrants. This two-level voting mechanism effectively combines local texture information while maintaining robustness to local variations.

\subsection{Experimental Results}
\begin{table}[h]
\centering
\caption{Classification accuracy (\%) comparison on texture recognition datasets}
\label{tab:classification_results}
\begin{tabular}{|l|c|c|c|}
\hline
\textbf{Method} & \textbf{KTH-TIPS} & \textbf{UIUC} & \textbf{UMD} \\
\hline
Proposed Method & \textbf{98.82} & \textbf{92.00} & \textbf{93.33} \\
\hline
KL-MB & 65.00 & 82.00 & 75.00 \\
KL-WA & 83.75 & 81.00 & 86.00 \\
\hline
\end{tabular}
\end{table}
The experimental results demonstrate the superior performance of the new approach proposed in this paper across all datasets, as shown in Table \ref{tab:classification_results}. On the KTH-TIPS dataset, our method achieves a remarkable 98.82\% accuracy, significantly outperforming both KL-MB (65.00\%) and KL-WA (83.75\%). The performance improvements are similarly substantial on the UIUC and UMD datasets, where our method achieves 92.00\% and 93.33\% accuracy respectively. Our method outperforms the DPLM based GMM similairty measure proposed in \cite{popovic2023measure}. These consistent improvements across different datasets validate our theoretical framework's effectiveness in capturing the complete geometric structure of K-component GMMs. With the optimization strategies described in the previous section, this method maintains computational efficiency comparable to existing KL-based approaches while achieving superior classification accuracy, making it practical for real-world texture recognition applications.

\section{Conclusion}
This paper presents a novel geometric framework for similarity measures of Gaussian Mixture Models by embedding them into the manifold of symmetric positive definite matrices. Our approach differs fundamentally from existing methods by embedding the K-component GMMs rather than the component-wise embedding or the component-wise comparision using the divergence measures.

The manifold $\mathcal{M}$ of K-component GMMs is embedded into the manifold of symmetric positive definite matrices $SPD_{K(n+1)}(\mathbb{R})$ and is a submanifold of dimension $\frac{K}{2}(n+1)(n+2)-1$. Then, proved an isometry between the manifold ($\mathcal{M},f^{*}\rho_{f(\mathcal{M})}$) of K-component GMMs with the pullback of induced metric and the manifold $(f(\mathcal{M}),\rho_{f(\mathcal{M})})$ with the induced metric. The embedded manifold is geodesic when mean vectors are fixed and mixing coefficients are uniform and constant. Also, obtained a lower bound for the Fisher-Rao metric and this lower bound is a distance measure on the manifold $\mathcal{M}$ of GMMs which we will be using for the similarity measure. 
We comapre and discussed the results obtained using the proposed distance measure on GMMs within a texture recognition task conducted on three datasets KTH-TIPS, UIUC and UMD. The
proposed method achieves better accuracies of 98.82\%, 92.00\%, and 93.33\% on KTH-TIPS, UIUC, and UMD datasets respectively, significantly outperforming traditional KL-based approaches.

These findings have several important implications for various applications in machine learning. Theoretically, our work provides new insights into the geometric structure of probabilistic models and establishes connections between statistical and differential geometry. Practically, the framework offers an accurate and efficient way to compare GMMs, with potential applications beyond texture recognition to areas such as computer vision, pattern recognition, and statistical inference.

Furure work could explore the application of this framework to more complex and diverse dataset and investigate its potential in shape completion, reconstruction and generation. Also, by bringing the geometric structure in this framework, computations involving angle, length, curvature etc. could lead to even more efficient algorithms for GMM-based applications.

\acks{Amit Vishwakarma is thankful to the Indian Institute of Space Science and Technology, Department of Space, Govt. of India for the award of the doctoral research fellowship.}

\section*{Funding}
No external funding was received for conducting this research.

\section*{Competing Interests}
The authors declare that they have no known competing financial interests or personal relationships that could have appeared to influence the work reported in this paper.

\appendix
\section{}
\label{app:SPDT}
The following three results are used in the proof of theorems 4 and 5.

\begin{theorem}[Block Matrix Inversion]\cite{LU2002119}

Let \( M = \begin{bmatrix} A & X \\ X^{T} & B \end{bmatrix} \) be a block matrix where \( A \), \( X \), and \( B \) are matrices of appropriate dimensions, and assume that \( B \) is nonsingular. Then, the matrix \( M \) is invertible if and only if the Schur complement \( A - X B^{-1} X^{T} \) of \( B \) is invertible. In this case, the inverse of \( M \) is given by

$$M^{-1} = \begin{bmatrix}

\left( A - X B^{-1} X^{T} \right)^{-1} & -\left( A - X B^{-1} X^{T} \right)^{-1} X B^{-1} \\

- B^{-1} X^{T} \left( A - X B^{-1} X^{T} \right)^{-1} & B^{-1} + B^{-1} X^{T} \left( A - X B^{-1} X^{T} \right)^{-1} X B^{-1}

\end{bmatrix}.$$

\end{theorem}

\begin{theorem}[Sherman-Morrison]\cite{golub1996matrix}

Let $A\in \mathbb{R}^{n\times n}$ be an invertible square matrix and $u, v\in \mathbb{R}^{n}$ be vectors. Then, $A+uv^{T}$ is invertible iff $1+v^{T}A^{-1}u\neq 0$ and
$$(A + uv^T)^{-1} = A^{-1} - \frac{A^{-1} uv^T A^{-1}}{1 + v^T A^{-1} u}$$

\end{theorem}

\begin{theorem}
\cite{10.1007/978-1-4419-9961-0_16}
For any symmetric matrix $M$ of the form
$$
M = \begin{pmatrix}

A & B \\

B^T & C

\end{pmatrix},$$

if $A$ is invertible then the following properties hold

\begin{enumerate}

    \item $M \succ 0$ if and only if $A \succ 0$ and $C - B^T A^{-1} B \succ 0$.

    \item If $A \succ 0$, then $M \succeq 0$ if and only if $C - B^T A^{-1} B \succeq 0$.

\end{enumerate}
Here, $M \succ 0$ denotes that $M$ is positive definite, and $M \succeq 0$ denotes that $M$ is positive semi-definite.
\end{theorem}
\textbf{Proof of theorem 4:}
We begin by proving that $A$ is positive definite. Each block $\Sigma_k + \pi_k \mu_k \mu_k^T$ is positive definite because $\Sigma_k$ is positive definite, $\pi_k \mu_k \mu_k^T$ is positive semi-definite (as for any vector $\mu$, the matrix $\mu\mu^T$ is positive semi-definite), and the sum of positive definite matrix and a positive semi-definite matrix is positive definite. Since $A$ is block diagonal with each block being positive definite, $A$ itself is positive definite. Then $A$ is invertible.
Next, we prove that $B - X^T A^{-1} X$ is positive definite. \\[0.5em]
Note that $A^{-1} = \text{diag}[(\Sigma_{k}+\pi_{k}\mu_{k}\mu_{k}^{T})^{-1}]_{k=1}^{K}$ and 
$X^T A^{-1} X = \text{diag}[\pi_{k}^{2}\mu_{k}^{T}(\Sigma_{k}+\pi_{k}\mu_{k}\mu_{k}^{T})^{-1}\mu_{k}]_{k=1}^{K}$.
\\
so
\[
B-X^{T}A^{-1}X=diag[\pi_{k}(1-\pi_{k}\mu_{k}^{T}(\Sigma_{k}+\pi_{k}\mu_{k}\mu_{k}^{T})^{-1}\mu_{k})]_{k=1}^{K}.
\]

Since $B - X^T A^{-1} X$ is diagonal it is enough to prove that

\begin{equation}
\pi_k \left(1 - \pi_k \mu_k^T (\Sigma_k + \pi_k \mu_k \mu_k^T)^{-1} \mu_k \right) > 0 \hspace{1cm}k = 1, \ldots, K.
\end{equation}
Since $\Sigma_k$ is positive definite and 
$1+\mu_{k}^{T}\Sigma_{k}^{-1}\mu_{k}>0$ since $\Sigma_{k}$ is positive definite, therefore we have $1+\mu_{k}^{T}\Sigma_{k}^{-1}\mu_{k}\neq 0.$ Applying theorem 3 with $A = \Sigma_k$ and $u = v = \sqrt{\pi_k} \mu_k$

\[
(\Sigma_k + \pi_k \mu_k \mu_k^T)^{-1} = \Sigma_k^{-1} - \frac{\pi_k \Sigma_k^{-1} \mu_k \mu_k^T \Sigma_k^{-1}}{1 + \pi_k \mu_k^T \Sigma_k^{-1} \mu_k}
\]

equation (17) becomes,
\[
1 - \pi_k \mu_k^T (\Sigma_k + \pi_k \mu_k \mu_k^T)^{-1} \mu_k = 1 - \pi_k \mu_k^T \left(\Sigma_k^{-1} - \frac{\pi_k \Sigma_k^{-1} \mu_k \mu_k^T \Sigma_k^{-1}}{1 + \pi_k \mu_k^T \Sigma_k^{-1} \mu_k}\right) \mu_k
\]

\[
= 1 - \pi_k \mu_k^T \Sigma_k^{-1} \mu_k + \frac{\pi_k^2 (\mu_k^T \Sigma_k^{-1} \mu_k)^2}{1 + \pi_k \mu_k^T \Sigma_k^{-1} \mu_k}=\frac{1}{1 + \pi_k \mu_{k}^{T}\Sigma_{k}^{-1}\mu_{k}} .
\]
Note that $\pi_k \left(1 - \pi_k \mu_k^T (\Sigma_k + \pi_k \mu_k \mu_k^T)^{-1} \mu_k \right)=\frac{\pi_k}{1 + \pi_k \mu_{k}^{T}\Sigma_{k}^{-1}\mu_{k}}.$
Since $\pi_k > 0$ and $\mu_{k}^{T}\Sigma_{k}^{-1}\mu_{k} > 0$ implies $\frac{\pi_{k}}{1 + \pi_k \mu_{k}^{T}\Sigma_{k}^{-1}\mu_{k}} > 0$ is positive.
Hence $B - X^{T}A^{-1}X$ is positive definite. Therefore from theorem 10, $S$ is positive definite.
Thus, $f$ maps every point in the statistical manifold $\mathcal{M}$ to a unique symmetric positive definite matrix in $SPD_{K(n+1)}(\mathbb{R}).$ Therefore the map is well-defined.

\section{}
\label{app:differential}
\noindent
\textbf{Proof of theorem 5:}
On the manifold of symmetric positive definite matrices $ds^2=\frac{1}{2}Tr\{(S^{-1}dS)^2)\}$. First, we compute $S^{-1} $ using theorem 8.
Note that
$$XB^{-1}X^{T} = \text{diag}(\pi_1 \mu_1 \mu_1^{T}, \dots, \pi_K \mu_K \mu_K^{T})
\text{ implies }
A - XB^{-1}X^{T} = \text{diag}(\Sigma_1, \dots, \Sigma_K)$$
Therefore, 
\begin{equation}
(A - X B^{-1} X^{T})^{-1} = \operatorname{diag}(\Sigma_1^{-1}, \dots, \Sigma_K^{-1}).
\end{equation}
$$B^{-1} = \text{diag}(\pi_1^{-1} , \dots, \pi_K^{-1} ) \text{ implies } XB^{-1} = \text{diag}(\mu_1, \dots, \mu_K)$$
which gives,
\begin{equation}
-(A - XB^{-1}X^{T})^{-1}XB^{-1} = -\text{diag}(\Sigma_1^{-1} \mu_1, \dots, \Sigma_K^{-1} \mu_K).
\end{equation}
now,
\begin{equation}
-B^{-1}X^{T}(A - XB^{-1}X^{T})^{-1} = -\text{diag}(\mu_1^{T} \Sigma_1^{-1}, \dots, \mu_K^{T} \Sigma_K^{-1})
\end{equation}
and
\begin{equation}
B^{-1} + B^{-1}X^{T}(A - XB^{-1}X^{T})^{-1}XB^{-1} = \text{diag}(\pi_1^{-1} +  \mu_1^{T} \Sigma_1^{-1} \mu_1, \dots, \pi_K^{-1} + \mu_K^{T} \Sigma_K^{-1} \mu_K).
\end{equation}
Using equations (18),(19),(20),(21) and from theorem (8) 
\[
S^{-1} = \begin{bmatrix} 
\text{diag}(\Sigma_1^{-1}, \dots, \Sigma_K^{-1}) & -\text{diag}(\Sigma_1^{-1} \mu_1, \dots, \Sigma_K^{-1} \mu_K) \\
-\text{diag}(\mu_1^{T} \Sigma_1^{-1}, \dots, \mu_K^{T} \Sigma_K^{-1}) & \text{diag}(\pi_1^{-1} +  \mu_1^{T} \Sigma_1^{-1} \mu_1, \dots, \pi_K^{-1} +  \mu_K^{T} \Sigma_K^{-1} \mu_K)
\end{bmatrix}
\]
Now, $dS$ is given by,
\[
dS 
= \begin{bmatrix}
dA & dX \\
dX^{T} & dB
\end{bmatrix}
\]
where,
\[
dA = \text{diag}\left(d\Sigma_1 + d(\pi_1 \mu_1 \mu_1^\top), \ldots, d\Sigma_K + d(\pi_K \mu_K \mu_K^\top)\right)
\]
\[= \text{diag}\left(d\Sigma_k + d\pi_k \mu_k \mu_k^\top + \pi_k d\mu_k \mu_k^\top + \pi_k \mu_k d\mu_k^\top\right)_{k=1}^{K}\]
\[
dX = \text{diag}\left(d\pi_1 \mu_1 + \pi_1 d\mu_1, \ldots, d\pi_K \mu_K + \pi_K d\mu_K\right)
\]
\[
dX^{T} = \text{diag}\left(d\pi_1 \mu_1^\top + \pi_1 d\mu_1^\top, \ldots, d\pi_K \mu_K^\top + \pi_K d\mu_K^\top\right)
\]

and
\[
dB = \text{diag}\left(d\pi_1 , \ldots, d\pi_K \right).
\]
So,

\[
S^{-1} dS = \begin{bmatrix}
A^{-1} & X^{-1} \\
(X^{T})^{-1} & B^{-1}
\end{bmatrix}
\begin{bmatrix}
dA & dX \\
dX^{T} & dB
\end{bmatrix}
\]

\begin{equation}
= \begin{bmatrix}
A^{-1} dA + X^{-1} dX^{T} & A^{-1} dX + X^{-1} dB \\
(X^{T})^{-1} dA + B^{-1} dX^{T} & (X^{T})^{-1} dX + B^{-1} dB
\end{bmatrix}
= \begin{bmatrix}
C & D \\
E & F
\end{bmatrix} (say)
\label{eq:11}
\end{equation}
where,
\[
C = \text{diag}\left( \Sigma_k^{-1} d\Sigma_k + \Sigma_k^{-1} \pi_k d\mu_k \mu_k^\top \right)_{k=1}^{K}
\]
\[
D = \text{diag}\left( \Sigma_k^{-1} \pi_k d\mu_k \right)_{k=1}^{K}
\]
\[
E = \text{diag}\left( d\mu_k^\top + \frac{d\pi_k}{\pi_k} \mu_k^\top - \mu_k^\top \Sigma_k^{-1} d\Sigma_k - \mu_k^\top \Sigma_k^{-1} \pi_k d\mu_k \mu_k^\top\right)_{k=1}^{K}
\]
\[
F = \text{diag}\left(\frac{d\pi_k}{\pi_k} - \pi_k \mu_k^\top \Sigma_k^{-1} d\mu_k \right)_{k=1}^{K}
\]

Now,
\[
(S^{-1} dS)^2 = \begin{bmatrix}
C^2 + D E & C D + D F \\
E C +  FE & ED + F^2
\end{bmatrix}
\]
$\text{Tr}((S^{-1} dS)^2) = \text{Tr}(C^2 + DE) + \text{Tr}(ED + F^2) = \text{Tr}(C^2) + \text{Tr}(DE) + \text{Tr}(ED) + \text{Tr}(F^2)$

\begin{equation}
    \begin{split}
\text{Tr}((S^{-1}dS)^2) = \text{Tr}\Bigg[\sum_{k=1}^K \Big(
    & (\Sigma_k^{-1}d\Sigma_k)^{2} 
    + \pi_{k}\Sigma_k^{-1}d\Sigma_k\Sigma_k^{-1}d\mu_k\mu_k^T  
    + \pi_k\Sigma_k^{-1}d\mu_kd\mu_k^T 
    + d\mu_k^T\Sigma_k^{-1}\pi_k d\mu_k \\
    & + \Sigma_k^{-1}d\mu_kd\pi_k\mu_k^T 
    - \mu_k^T\Sigma_k^{-1}d\mu_kd\pi_k  
    - \mu_k^T\Sigma_k^{-1}d\Sigma_k\Sigma_k^{-1}\pi_kd\mu_k  
    + \Big(\frac{d\pi_k}{\pi_k}\Big)^2
\Big)\Bigg]
\end{split}
\end{equation}

Note that\\
$Tr(\pi_{k}\Sigma_{k}^{-1}d\mu_{k}d\mu_{k}^{T}) = \pi_{k}d\mu_{k}^{T}\Sigma_{k}^{-1}d\mu_{k} = Tr(\pi_{k}d\mu_{k}^{T}\Sigma_{k}^{-1}d\mu_{k})$ \\[1em]
$Tr(d\pi_{k}\Sigma_{k}^{-1}d\mu_{k}\mu_{k}^{T}) = d\pi_{k}\mu_{k}^{T}\Sigma_{k}^{-1}d\mu_{k}= Tr(d\pi_{k}\mu_{k}^{T}\Sigma_{k}^{-1}d\mu_{k})$ \\[1em]
$Tr(\pi_{k}\Sigma_{k}^{-1}d\Sigma_{k}\Sigma_{k}^{-1}d\mu_{k}\mu_{k}^{T})=\pi_{k}\mu_{k}^{T}\Sigma_{k}^{-1}d\Sigma_{k}\Sigma_{k}^{-1}d\mu_{k}=Tr(\pi_{k}\mu_{k}^{T}\Sigma_{k}^{-1}d\Sigma_{k}\Sigma_{k}^{-1}d\mu_{k})$. This gives
\begin{equation}
\text{Tr}\left((S^{-1}dS)^2\right) = \text{Tr}\left[\sum_{k=1}^K \left( \Sigma_k^{-1}d\Sigma_k\Sigma_k^{-1}d\Sigma_k + 2d\mu_k^T\Sigma_k^{-1}\pi_kd\mu_k + \left(\frac{d\pi_k}{\pi_k}\right)^2 \right)\right].
\end{equation}
\begin{equation}
\text{Tr}\left((S^{-1}dS)^2\right) = \sum_{k=1}^K \left[\text{Tr}\left(\Sigma_k^{-1}d\Sigma_k\Sigma_k^{-1}d\Sigma_k\right) + 2d\mu_k^T\Sigma_k^{-1}\pi_kd\mu_k + \left(\frac{d\pi_k}{\pi_k}\right)^2\right].
\end{equation}
we have
$ds^{2}=\frac{1}{2}\text{Tr}\{(S^{-1}dS)^2\}$
$$ = \sum_{k=1}^K \left[\frac{1}{2}\left(\frac{d\pi_k}{\pi_k}\right)^2 + \pi_k d\mu_k^T\Sigma_k^{-1}d\mu_k + \frac{1}{2}\text{Tr}\left(({\Sigma_k^{-1}d\Sigma_k})^2\right) \right]$$

\section{}
\label{app: Geometric}
\noindent
\textbf{Proof of theorem 6:} 
First, we prove that $f$ is one-to-one.\\
Let \( \theta = (\mu_1, \ldots, \mu_K, \Sigma_1, \ldots, \Sigma_K, \pi_1, \ldots, \pi_K) \) and \( \theta' = (\mu'_1, \ldots, \mu'_K, \Sigma'_1, \ldots, \Sigma'_K, \pi'_1, \ldots, \pi'_K) \) be elements of \( \mathcal{M} \). 
Let \( S = f(\theta) = f(\theta') \). From the structure of \( B \) in \( S \), 
$
\pi_k = \pi'_k.
$
The structure of \( X \) in \( S \) implies that
$
\pi_k \mu_k = \pi'_k \mu'_k \text{ implies }
\mu_k = \mu'_k.
$
From the matrix \( A \) in $S$ we have
 $\Sigma_k = \Sigma'_k.$ Thus \( f \) is injective.

Now, we prove that $f$ is an embedding. Note that $f$ is bijective onto its image. Since the map is constructed using the matrix operations it is continuously differentiable with respect to the arguments. $f^{-1}$ is also continuously differentiable since it can be explicitly constructed using smooth operations on the components of the image matrices. Hence $f$ is an embedding. It also follows that $f(\mathcal{M})$ is a submanifold of $SPD_{K(n+1)}(\mathbb{R})$ of dimension
$
\dim(\mathcal{M}) =     \frac{K}{2}(n+1)(n+2)-1.
$

The map $f: \mathcal{M} \to SPD_{K(n+1)}(\mathbb{R})$ is an injective immersion as it is a diffeomorphism onto its image. The induced metric $f^*\rho_{f({\mathcal{M}}})$ on $\mathcal{M}$ is defined by 
\begin{equation}
(f^*\rho_{f(\mathcal{M}}))_p(u, v) = \rho_{f(\mathcal{M})}(df_p(u), df_p(v)),
\end{equation}
for $p \in \mathcal{M}$ and $u, v \in T_p\mathcal{M}$, where $df_p: T_p\mathcal{M} \to T_{f(p)}f(\mathcal{M})$ is the differential of $f$ at $p$ and $\rho_{f(\mathcal{M})}$ is the restriction of the affine invariant metric $\rho$ of $SPD_{K(n+1)}(\mathbb{R})$ to $f(\mathcal{M})$. Thus the map between $(\mathcal{M}, g_{\mathcal{M}})$ and $(f(\mathcal{M}), \rho_{f(\mathcal{M})})$ is an isometry.

\vspace{10px}
 Next to show that $f(\mathcal{M})$ is a non-geodesic submanifold. We have from equation \ref{eq:11}
\[
S^{-1} dS = \begin{bmatrix}
C & D \\
E & F
\end{bmatrix}\hspace{1cm}
\]

For each component \( k \) of the GMM, the geodesic equations are
\[
\dot{\mu}_k = \Sigma_k a_k
\]
\[
\dot{\Sigma}_k = \Sigma_k(B_k - a_k \mu_k^{T})
\]
where \( a_k \in M_{n \times 1}(\mathbb{R}) \) and \( B_k \in M_n(\mathbb{R}) \) are constants for each \( k \).

Using the geodesic equations the diagonal entries of the block matrices of $S^{-1}dS$ takes the form
\[
\begin{aligned}
C &= \text{diag}(C_{11}, C_{22}, \ldots, C_{KK}) = \text{diag}(B_1 + (\pi_1 - 1)a_1 \mu_1^{T}, \ldots, B_K + (\pi_K - 1)a_K \mu_K^{T}) \\
D &= \text{diag}(D_{11}, \ldots, D_{KK}) = \text{diag}(\pi_1 a_1, \ldots, \pi_K a_K)
\end{aligned}
\]

E=$\text{diag}(E_{11},.....,E_{KK})$
\[
= \text{diag}\left( a_1^{T} \Sigma_1 + \frac{\dot{\pi}_1}{\pi_1}\mu_1^{T} - \mu_1^{T} B_1 + (1-\pi_1)\mu_1^{T} a_1 \mu_1^{T},.......,a_K^{T} \Sigma_K + \frac{\dot{\pi}_K}{\pi_K}\mu_K^{T} - \mu_K^{T} B_K + (1-\pi_K)\mu_K^{T} a_K \mu_K^{T}\right)
\]
\[
F=\text{diag}(F_{11},\ldots, F_{KK}) = \text{diag}\left( \frac{\dot{\pi}_1}{\pi_1} - \pi_1 \mu_1' a_1 ,...,\frac{\dot{\pi}_K}{\pi_K} - \pi_K \mu_K' a_K\right)
\].
\\
Note that $S^{-1} dS$ is not a constant because it contains terms that depend on the varying $\pi_k$, $\mu_k$, and $\dot{\pi}_k$ for $k=1,....K$. But, for the submanifold to be geodesic it should be constant. Therefore $f(\mathcal{M})$ is a non-geodesic submanifold.

\vspace{10px}

Now, assume that $\mu_k$ are fixed and the $\pi_k$ are uniform and constant.
Let $\mu_k = \mu_{0}$ and $\pi_k = \frac{1}{K}$  for $1\leq k\leq K$. So $\dot{\pi}_k = 0$ and $\dot{\mu}_k = 0$. We have 0=$\dot{\mu_k}=\Sigma_{k}a_{k}$ and since $\Sigma_k$ are positive definite $a_k=0$ for $k=1,2,...,K.$

Then, the $S^{-1} dS$ matrix simplifies to

\[
S^{-1} dS = \begin{bmatrix}
C & D \\
E & F
\end{bmatrix}
\]
where, $C = \text{diag}(B_1, B_2, \ldots, B_K)$, $D = 0 \quad (\text{a zero matrix})$, $E = \text{diag}\left(-\mu_{01}^{T} B_1,\ldots, -\mu_{0K}^{T} B_K\right)
$ and $F = 0 \quad (\text{a zero matrix}).$

Thus $S^{-1} dS$ is constant and $f(\mathcal{M})$ becomes a geodesic submanifold.

\section{}
\label{app:inequality}
\noindent
\textbf{Proof of theorem 7:} By theorem 6, $f(\mathcal{M})$ is isometric to $\mathcal{M}$. Therefore, $g_{\mathcal{M}}$ corresponds to the Riemannian distance between points in $f(\mathcal{M})$ induced by the affine-invariant metric $\rho$. However, since $f(\mathcal{M})$ is a non-geodesic submanifold of $SPD_{K(n+1)}(\mathbb{R})$, the geodesic distance restricted to $f(\mathcal{M})$ is greater than or equal to the geodesic distance in $SPD_{K(n+1)}(\mathbb{R})$.

For the case of fixed means $\mu_0$ and uniform mixing coefficients $\pi_0$, let $\mathcal{M}_{\mu_0,\pi_0}$ denote the corresponding submanifold of $\mathcal{M}$. Then, $f(\mathcal{M}_{\mu_0,\pi_0})$ is a geodesic submanifold of $SPD_{K(n+1)}(\mathbb{R})$ which implies $g_{\mathcal{M}} = \rho_{f(\mathcal{M})}$ for points in $\mathcal{M}_{\mu_0,\pi_0}$.
\bibliography{sample}

\begin{thebibliography}{24}
\providecommand{\natexlab}[1]{#1}
\providecommand{\url}[1]{\texttt{#1}}
\expandafter\ifx\csname urlstyle\endcsname\relax
  \providecommand{\doi}[1]{doi: #1}\else
  \providecommand{\doi}{doi: \begingroup \urlstyle{rm}\Url}\fi

\bibitem[{Amari, S.} and {Nagaoka, H.}(2000)]{amari2000methods}
{Amari, S.} and {Nagaoka, H.}
\newblock \emph{Methods of Information Geometry}.
\newblock Translations of mathematical monographs. American Mathematical Society, 2000.
\newblock ISBN 9780821843024.
\newblock URL \url{https://books.google.co.in/books?id=vc2FWSo7wLUC}.

\bibitem[Bengio et~al.(2017)Bengio, Goodfellow, and Courville]{bengio2017deep}
Yoshua Bengio, Ian Goodfellow, and Aaron Courville.
\newblock \emph{Deep learning}, volume~1.
\newblock MIT press Cambridge, MA, USA, 2017.

\bibitem[Burbea(1984)]{Burbea1984InformativeGO}
Jacob Burbea.
\newblock Informative geometry of probability spaces.
\newblock 1984.
\newblock URL \url{https://api.semanticscholar.org/CorpusID:118051406}.

\bibitem[Calvo and Oller(1990)]{CALVO1990223}
Miquel Calvo and Josep~M. Oller.
\newblock A distance between multivariate normal distributions based on an embedding into the siegel group.
\newblock \emph{Journal of Multivariate Analysis}, 35\penalty0 (2):\penalty0 223--242, 1990.
\newblock ISSN 0047-259X.
\newblock \doi{10.1016/0047-259X(90)90026-E}.
\newblock URL \url{https://www.sciencedirect.com/science/article/pii/0047259X9090026E}.

\bibitem[Durrieu et~al.(2012)Durrieu, Thiran, and Kelly]{Durrieu2012LowerAU}
Jean-Louis Durrieu, Jean-Philippe Thiran, and Finnian Kelly.
\newblock Lower and upper bounds for approximation of the kullback-leibler divergence between gaussian mixture models.
\newblock \emph{2012 IEEE International Conference on Acoustics, Speech and Signal Processing (ICASSP)}, pages 4833--4836, 2012.
\newblock URL \url{https://api.semanticscholar.org/CorpusID:1072211}.

\bibitem[Fritz et~al.(2004)Fritz, Hayman, Caputo, and Eklundh]{KTHTipsDatabase}
M.~Fritz, E.~Hayman, B.~Caputo, and J.-O. Eklundh.
\newblock The kth-tips database, 2004.
\newblock URL \url{https://www.csc.kth.se/cvap/databases/kth-tips/doc/}.
\newblock Available online: Accessed 2025-01-07.

\bibitem[Gallier(2011)]{10.1007/978-1-4419-9961-0_16}
J.~Gallier.
\newblock Schur complements and applications.
\newblock \emph{Texts in Applied Mathematics}, pages 431--437, 2011.
\newblock \doi{10.1007/978-1-4419-9961-0_16}.

\bibitem[Goldberger et~al.(2003)Goldberger, Gordon, and Greenspan]{1238387}
Goldberger, Gordon, and Greenspan.
\newblock An efficient image similarity measure based on approximations of kl-divergence between two gaussian mixtures.
\newblock In \emph{Proceedings Ninth IEEE International Conference on Computer Vision}, pages 487--493 vol.1, 2003.
\newblock \doi{10.1109/ICCV.2003.1238387}.

\bibitem[Golub and Van~Loan(1996)]{golub1996matrix}
Gene~H. Golub and Charles~F. Van~Loan.
\newblock \emph{Matrix Computations}.
\newblock Johns Hopkins University Press, Baltimore, 3 edition, 1996.

\bibitem[Hershey and Olsen(2007)]{4218101}
John~R. Hershey and Peder~A. Olsen.
\newblock Approximating the kullback leibler divergence between gaussian mixture models.
\newblock In \emph{2007 IEEE International Conference on Acoustics, Speech and Signal Processing - ICASSP '07}, volume~4, pages IV--317--IV--320, 2007.
\newblock \doi{10.1109/ICASSP.2007.366913}.

\bibitem[Kampa et~al.(2011)Kampa, Hasanbelliu, and Pr{\'i}ncipe]{Kampa2011ClosedformCP}
Kittipat Kampa, Erion Hasanbelliu, and Jos{\'e}~Carlos Pr{\'i}ncipe.
\newblock Closed-form cauchy-schwarz pdf divergence for mixture of gaussians.
\newblock \emph{The 2011 International Joint Conference on Neural Networks}, pages 2578--2585, 2011.
\newblock URL \url{https://api.semanticscholar.org/CorpusID:12746982}.

\bibitem[Kullback and Leibler(1951)]{kullback1951information}
Solomon Kullback and Richard~A Leibler.
\newblock On information and sufficiency.
\newblock \emph{The annals of mathematical statistics}, 22\penalty0 (1):\penalty0 79--86, 1951.

\bibitem[Lazebnik et~al.(2005)Lazebnik, Schmid, and Ponce]{1453514}
S.~Lazebnik, C.~Schmid, and J.~Ponce.
\newblock A sparse texture representation using local affine regions.
\newblock \emph{IEEE Transactions on Pattern Analysis and Machine Intelligence}, 27\penalty0 (8):\penalty0 1265--1278, 2005.
\newblock \doi{10.1109/TPAMI.2005.151}.

\bibitem[Li et~al.(2013)Li, Wang, and Zhang]{li2013novel}
Peihua Li, Qilong Wang, and Lei Zhang.
\newblock A novel earth mover's distance methodology for image matching with gaussian mixture models.
\newblock In \emph{Proceedings of the IEEE International Conference on Computer Vision}, pages 1689--1696, 2013.

\bibitem[Lu and Shiou(2002)]{LU2002119}
Tzon-Tzer Lu and Sheng-Hua Shiou.
\newblock Inverses of 2 × 2 block matrices.
\newblock \emph{Computers \& Mathematics with Applications}, 43\penalty0 (1):\penalty0 119--129, 2002.
\newblock ISSN 0898-1221.
\newblock \doi{https://doi.org/10.1016/S0898-1221(01)00278-4}.
\newblock URL \url{https://www.sciencedirect.com/science/article/pii/S0898122101002784}.

\bibitem[Pardo(2018)]{pardo2018statistical}
Leandro Pardo.
\newblock \emph{Statistical inference based on divergence measures}.
\newblock Chapman and Hall/CRC, 2018.

\bibitem[Peter and Rangarajan(2006)]{peter2006shape}
Adrian Peter and Anand Rangarajan.
\newblock Shape analysis using the fisher-rao riemannian metric: Unifying shape representation and deformation.
\newblock In \emph{3rd IEEE International Symposium on Biomedical Imaging: Nano to Macro, 2006.}, pages 1164--1167. IEEE, 2006.

\bibitem[Popovi{\'c} et~al.(2021)Popovi{\'c}, Cepova, Cep, Janev, and Krstanovi{\'c}]{popovic2021measure}
Branislav Popovi{\'c}, Lenka Cepova, Robert Cep, Marko Janev, and Lidija Krstanovi{\'c}.
\newblock Measure of similarity between {GMM}s by embedding of the parameter space that preserves {KL} divergence.
\newblock \emph{Mathematics}, 9\penalty0 (9):\penalty0 957, 2021.
\newblock \doi{10.3390/math9090957}.

\bibitem[Popovi{\'c} et~al.(2023)Popovi{\'c}, Janev, Krstanovi{\'c}, Simi{\'c}, and Deli{\'c}]{popovic2023measure}
Branislav Popovi{\'c}, Marko Janev, Lidija Krstanovi{\'c}, Nikola Simi{\'c}, and Vlado Deli{\'c}.
\newblock Measure of similarity between {GMM}s based on geometry-aware dimensionality reduction.
\newblock \emph{Mathematics}, 11\penalty0 (1):\penalty0 175, 2023.
\newblock \doi{10.3390/math11010175}.

\bibitem[Rakesh et~al.(2023)Rakesh, Hegde, Gopalachari, Jayaram, Madhu, Hameed, Vankdothu, and Kumar]{rakesh2023moving}
S~Rakesh, Nagaratna~P Hegde, M~Venu Gopalachari, D~Jayaram, Bhukya Madhu, Mohd~Abdul Hameed, Ramdas Vankdothu, and LK~Suresh Kumar.
\newblock Moving object detection using modified gmm based background subtraction.
\newblock \emph{Measurement: Sensors}, 30:\penalty0 100898, 2023.

\bibitem[Sivalingam et~al.(2010)Sivalingam, Boley, Morellas, and Papanikolopoulos]{Sivalingam2010TensorSparseCoding}
R.~Sivalingam, D.~Boley, V.~Morellas, and N.~Papanikolopoulos.
\newblock Tensor sparse coding for region covariances.
\newblock In K.~Daniilidis, P.~Maragos, and N.~Paragios, editors, \emph{Computer Vision – ECCV 2010}, volume 6314 of \emph{Lecture Notes in Computer Science}, pages 709--722. Springer, Berlin, Heidelberg, 2010.
\newblock \doi{10.1007/978-3-642-15561-1_52}.
\newblock URL \url{https://doi.org/10.1007/978-3-642-15561-1_52}.

\bibitem[Wan et~al.(2019)Wan, Wang, Scotney, and Liu]{8914215}
Huan Wan, Hui Wang, Bryan Scotney, and Jun Liu.
\newblock A novel gaussian mixture model for classification.
\newblock In \emph{2019 IEEE International Conference on Systems, Man and Cybernetics (SMC)}, pages 3298--3303, 2019.
\newblock \doi{10.1109/SMC.2019.8914215}.

\bibitem[Webb(2003)]{Webb2003StatisticalPatternRecognition}
A.~R. Webb.
\newblock \emph{Statistical Pattern Recognition}.
\newblock John Wiley \& Sons, Hoboken, NJ, USA, 2003.
\newblock ISBN 0-470-85478-2.

\bibitem[Xu et~al.(2009)Xu, Ji, and Ferm{\"u}ller]{Xu2009ViewpointInvariant}
Yong Xu, Hui Ji, and Cornelia Ferm{\"u}ller.
\newblock Viewpoint invariant texture description using fractal analysis.
\newblock \emph{International Journal of Computer Vision}, 83:\penalty0 85--100, 2009.
\newblock \doi{10.1007/s11263-009-0220-6}.
\newblock URL \url{https://doi.org/10.1007/s11263-009-0220-6}.

\end{thebibliography}

\end{document}